\documentclass[preprint]{elsarticle}
\usepackage[utf8]{inputenc}
\usepackage{amsfonts}
\usepackage{amsmath}
\usepackage{algpseudocode}
\usepackage{mathtools}
\usepackage{comment}
\usepackage{url}
\usepackage{xcolor}
\usepackage{booktabs}
\newtheorem{definition}{Definition}
\newcommand{\maximize}{%
	\mathopen{}\operatorname*{max}%
}
\newcommand{\minimize}{%
	\mathopen{}\operatorname*{min}%
}
\title{Bilevel optimization for the deployment of refueling stations for electric vehicles on road networks}
\author[1]{Ramón Piedra-de-la-Cuadra\corref{cor1}%
	\fnref{fn1}}
\ead{rpiedra@us.es}

\author[1]{Francisco A. Ortega\fnref{fn3}}
\ead{riejos@us.es}
\cortext[cor1]{Corresponding author}
\fntext[fn1]{rpiedra@us.es}

\fntext[fn3]{riejos@us.es}
\affiliation[1]{organization={Departamento de Matemática Aplicada I, Universidad de Sevilla},
	addressline={Avda. Reina Mercedes, s/n},
	postcode={41012},
	city={Seville},
	country={Spain}}

\date{\today}

\begin{document}

	\begin{abstract}
		This work consists of a procedure to optimally select, among a group of candidate sites where gas stations were already located, a sufficient number of charging points in order to guarantee that an electric vehicle can make its journey without a problem of energy autonomy and that each selected charging station has another one that serves as useful support in case of failure (reinforced coverage service). For this purpose, we propose a bilevel model that, in a former level, minimizes the number of refueling points necessary to guarantee a reinforced service coverage for all users who transit from their origin to destination and, as a second level, maximize the volume of demand that can be satisfied subject to budgetary restrictions. With the first of the objectives we are addressing the typical {requirement} of the administration, which consists of guaranteeing the viability of the solutions, and the second of the objectives is a criterion typically used by the private sector initiative, compatible with the profit maximization.
	\end{abstract}
\begin{keyword}Bilevel optimization \sep Electric vehicles \sep Conditional coverage problem \sep Knapsack problem
\end{keyword}
\maketitle

	\section{Introduction}

\begin{figure}[h!]
		\centering
		\includegraphics[scale=0.075]{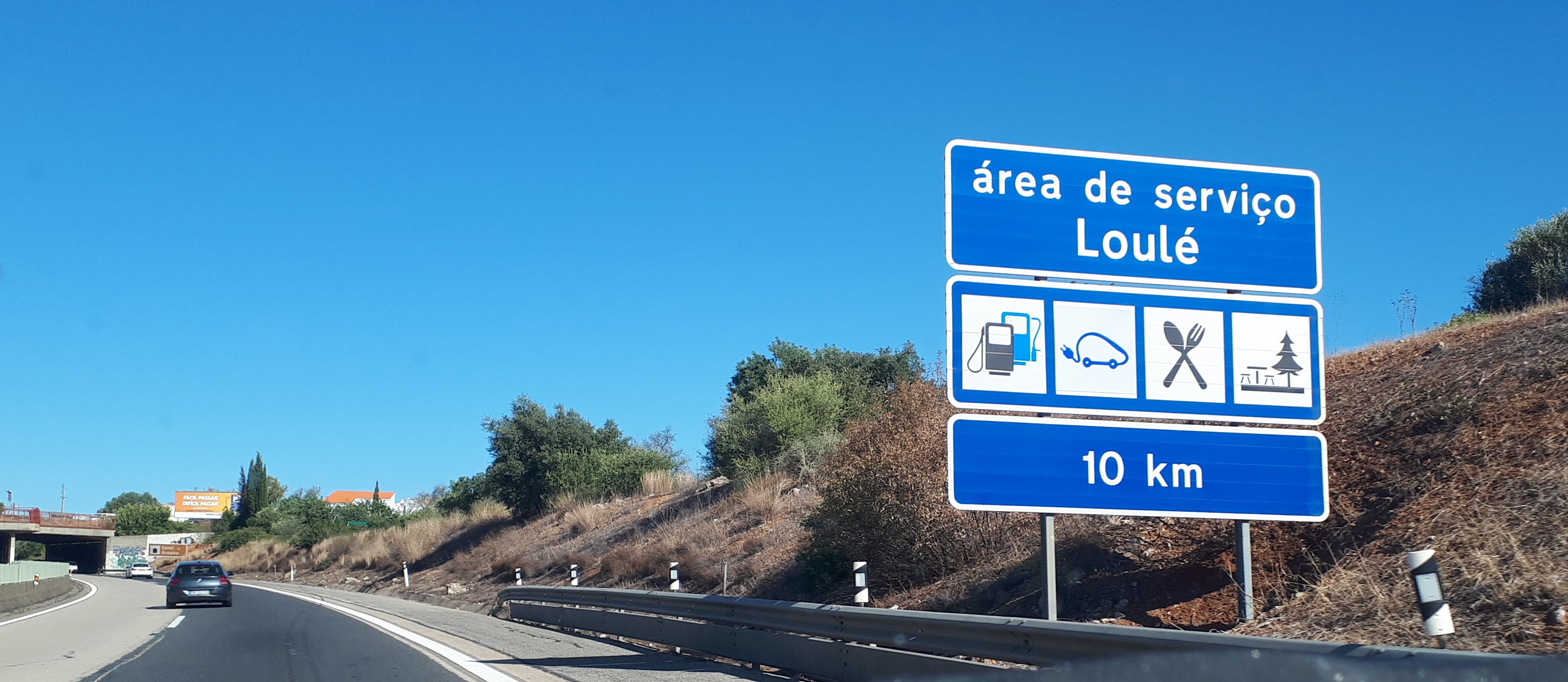}
		\caption{Example {of} adapted gas station}
		\label{Fig1}
	\end{figure}

	According to the Climate Plan decided at the 26th session of the Conference of the Parties (COP26) to the United Nations Framework Convention on Climate Change (UNFCCC), many governments have taken multiple strategic policy initiatives in the energy and transport sector to steer their respective nations to the path of reducing total projected carbon emissions by one billion tons from now to 2030. The European Union (EU) is the world's third-largest emitter of greenhouse gases behind China and the USA, followed by India and Russia. The European Commission has set out the objective of leading the world in the transition to a carbon neutral economy and established a goal of net-zero economy-wide emissions by 2050. \\
	Transport was responsible for close to a quarter of CO2 emissions in the EU in 2019, of which 71.7\% came from road transport, according to the annual European Environment Agency report. To reduce CO2 emissions and achieve the climate neutrality of the European Green Deal, greenhouse gas emissions from transport would have to be reduced by 90\% by 2050, compared to 1990 levels. However, current projections put the decline in transport emissions by 2050 at just 22\%, well below current targets.
	Electric cars (EC) have a series of very important advantages over conventional internal combustion engine models linked to fuel savings, tax exemption in some countries, subsidies, maintenance costs, and zero emission of polluting gases to the atmosphere. For these reasons, 
	the recent market shift towards electric vehicles (EVs) in Europe has been impressive. In 2020, despite the contraction of overall car sales in Europe, EV registrations more than doubled to 1.4 million and reached 10\% of the market, while this number stood at 6\% in China and 2\% in the US \cite{IEA21}.\\
	 It has been widely studied that availability of charging infrastructure is a major concern for adoption of electric vehicles \cite{Nicholas20}. {The availability of infrastructure to serve the electric car market varies considerably between EU Member States. For example, in the Netherlands there are more than 32,000 charging points and there are more than 119,000 registered electric vehicles; while in Greece, only 40 charging points are available for just over 300 EVs \cite{Niestadt19}.}\\
	According to the Spanish Electrical Network, the average mobility in Spain is 40 {km} per day.
	Depending on the characteristics of the trip and the vehicle, the need to recharge the battery will be essential at some point during the itinerary.\\
	An electric car charging station is a place where vehicles equipped with a plug (whether they are 100\% electric or hybrid) obtain the energy they need to function, in the same way that propulsion vehicles do with diesel or gasoline. The types of fast and ultra-fast recharge for EVs resemble the current refueling characteristics of combustion vehicles.  The vehicle must remain connected during some minutes to obtain an 80\% charge. Service stations or roadside restaurants located on main interurban roads, where the user can stay for a short period of time, are ideal places to install this type of fueling. In addition, consumers' perception of electric vehicles is that they cannot cover the desired intercity distance without recharging \cite{Berkeley17}. 
	 A well spread public charging network has also shown increased propensity to encourage complete shift to EV in the household \cite{Lorentzen17}.\\
	{In the context of Spain, all the plans since 2000 to increase the use of EVs have failed, and the objectives have not been achieved \cite{Borge22}. 
	The causes are intricate and multifaceted, encompassing various barriers, such as the absence of a well-developed recharging infrastructure, a convoluted tax system, the limited operational range of EVs', costly recharging tariffs, and the scarcity of affordable vehicle options, among other factors \cite{Rosales22}.} That {is} why Spain's Recovery, Transformation and Resilience Plan, a goal has been set for 2023 of at least 100,000 charging points and 250,000 electric vehicles, as well as the development of the value chain, new business models and new dynamics that {favour} the progressive electrification of mobility, the reduction of emissions and the {fulfilment} of energy and climate objectives. The financing planned in Spain for this purpose will reach 140,000 million euros in transfers and credits over the next six years.\\
	In the case of electric mobility, it is planned to deploy fast or ultra-fast recharging infrastructure along corridors that make it possible to structure the entire territory, both on interurban roads of special national and regional relevance, as well as on those corridors connected with {neighbouring} countries.\\
	The industrial companies that manage gas stations are currently in a general process of adaptation with the aim of installing charging points for electric cars in their service stations. These energy companies develop strategic plans to adapt their infrastructures with the aim that electric vehicle drivers can have refueling points distributed throughout the territory for recharging energy, where these ones being not too far from each other to make it possible to circulate with supply guarantees throughout any corridor contained in the territory. Figure \ref{Fig1} shows an informative sign indicating the proximity of an EV refueling station along a road in the Algarve (Portugal).
	\\
	From a double public-private vision, we propose to select locations that coincide with certain existing traditional gas stations, where several electric vehicle recharging devices will be available, grouped to allow several users to be served simultaneously in a similar way to how demand is met of petroleum-derived fuel at conventional gas stations.\\
	This work develops a procedure to optimally select, among a group of candidate sites where gas stations are already located, enough charging points in such a way as to guarantee that any electric vehicle can make its journey without a problem of energy autonomy and that each selected charging station has another at least that serves as coverage in case of failure (reinforced service). For this purpose, we propose a bilevel model that minimizes the number of refuelling points necessary to guarantee a reinforced service coverage for all users who transit from their origin to destination inside a territorial zone and, as a second level, maximize the volume of served demand subject to budgetary restrictions. With the first of these objectives, we are meeting the usual requirement of the administration, which consists of guaranteeing the viability of the solutions, and the second of the objectives is a criterion typically used by the private sector initiative of profit maximization.\\
	The remainder of the paper is organized as follows. Literature background of the work is outlined in Section 2. Section 3 formulates the problem of {planning an optimal deployment of} charging stations as a bi-level mathematical program and proposes its solution algorithm. Section 4 is devoted to the presentation of our case study. Other perspective for the problem, {consisting if reformulating the model as a single-level optimization problem by exploiting the primal-dual optimality conditions,} is discussed in Section 5. Conclusions and directions for further research finally follow in Section 6.

	\section{Literature Review}

	Research on EV and its related infrastructure has gained momentum in the last decade. 
	{Ref.} \cite{Pagany19} presents an extensive overview as well as an in-depth review of the literature dealing with the location of charging stations (CS) for EV.\\
	Types of economic approaches used to locate charging stations can be classified in theoretical models and empirical applications.  Theoretical works have concentrated on deriving mathematical models to assess locations for EV infrastructure with simulated data whereas empirical models focus on real spaces with spatial characters often with georeferenced data and provide the results obtained from different mathematical locational models. EV charging infrastructure site selection is a multiple criteria decision making problem as it is determined by several factors often contradictory. One part of existing literature presents strategies for locating charging stations; in particular, which type of CS and where should be placed those installations within a selected area or along a road network are investigated. While some proposals attempt to identify those locations by calculating the spatial distribution of charging demand, other studies aim to find the best CS location with an optimization approach using traffic demand along road networks \cite{Ghosh22}.\\
	The location of electric charging stations is a much-discussed topic in the literature. {In} \cite{Lee22} {are review} the location models that have been employed in charging facilities for EVs and alternative energy-powered vehicles are reviewed: p-median problem, set covering problem, fixed charge location problem, and those based on demand of Origin-Destination trips. Moreover, contributions regarding the sizing problem of EV charging stations with different objective functions are also examined.\\
	{Ref.} \cite{Kuby05} introduced the flow refuelling-location model (FRLM) to maximize the impact of a given investment in refuelling infrastructure for alt-fuel vehicles. The model optimizes the location of a given number of refuelling facilities within a network with the goal of enabling the maximum number of trips by vehicles with a limited range. The most basic input {for} the FRLM is a set of origin–destination (O–D) pairs and the flow volumes between them. For each O–D pair, one must calculate the shortest path between origin and destination and then determine all the combinations of facilities that can refuel a round trip along that path. These combinations of refuelling stations, and the flow volume associated with each combination, are evaluated by the FRLM’s mixed integer programing formulation to determine optimal facility locations.\\
	The FRLM is an un-capacitated model; it implicitly assumes that a single facility can refuel an infinite amount of flow. This may not be a realistic assumption. To address this concern, {Ref.} \cite{Kuby07} introduces the capacitated flow refuelling location model (CFRLM) that limits the amount of flow that any facility can refuel. The objective consists of locating p stations on a network to maximize the refuelling of origin–destination flows. Due to the limited driving range of vehicles, network vertices do not constitute a finite dominating set. Such authors propose to add candidate sites along arcs using minimax and maximin methods. The maximin criterion is motivated by the idea of not wasting candidate sites by putting any too close together, while the minimax objective aims to avoid long arcs with no candidate sites. Nevertheless, none of the methods reaches to generate a finite dominating set.\\
	{Ref.} \cite{Upchurch09} extended the FRLM by considering that only limited number of vehicles can be refuelled by a station. To formulate this constraint, a variable was introduced to indicate the proportion of the traffic flow on each path being refuelled by each station combination.\\
	Another formulation of the FRLM was proposed by \cite{Capar12}, which does not require the pre-generation of all feasible station combinations. For this purpose, authors introduced two binary variables on each node along every path indicating whether a station exists at that node and whether a driver at that node can reach another station further down the path without running out of fuel.\\
	Literature contains another type of contributions for the FRLP which is based on the set covering problem.
	A generalized maximum covering model is proposed by \cite{Wen14} for the {flow refuelling-location problem (FRLP)} without using extra variables, {as happens in} \cite{Capar12}, and without pre-generating facility combinations, as in the other maximum covering models. A set of sub-paths is defined for each path, in such a way that if each of these sub-paths contains a replenishment station, the entire path flow is captured.\\
 {Ref.} \cite{Huang15} extend the set cover model by allowing shortest path deviations, where the deviation paths are exogenously determined, and fuel level is 
	tracked on every node.\\
	 {Ref.} \cite{Wang09} presented a model to capture all the traffic flow with least station location cost. In their model, a variable is defined at each node on each path indicating the remaining amount of fuel when vehicles reach that node, such that a trip can be refuelled if the remaining amount of fuel at each node along the path is non-negative.\\
	{Ref.} \cite{Wang10} consider that each node is associated with a population and, if a station is located in the node, the corresponding population is said to be covered. For this context, the authors extended the FRLP model to multiple objectives: to minimize the cost of locating the stations and also to maximize the coverage of the population.\\
	Special work has been done from the point of view of the company, maximizing the flow captured or minimizing the cost of the infrastructures. Optimization procedures employ both exact methods \cite{Asamer16, Li16, Wang16} as well as heuristic techniques \cite{Chen15,Hidalgo16,Salmon16,Sebastiani16}. They aim to find the optimal location for CS by minimizing total cost, reducing trip length, or waiting time.\\
	The locations and sizes of fast-charging stations in a transportation network should satisfy EV driving demands, while simultaneously ensuring the security operation constraints of power systems, e.g., distribution line current limits and nodal voltage limits. In addition, authors consider that an appropriate fast-charging station planning method should minimize the investment costs of both charging stations and corresponding power grid upgrades.\\
	In \cite{Motoaki19} only a location problem with the objective of providing a geographical coverage of the demand is considered. In this case, the objective is to have a maximum number of EVs with access to a potentially available station, and the charging station locations are uncorrelated to the charging station sizes.\\
	Bilevel modelling has also been used in the literature to study different objectives in the location of charging stations. {Ref.} \cite{Jing17} developed a bi-level model to maximize coverage of EV flows by deploying a given number of charging stations on a network with mixed conventional vehicles and EVs.\\
	{Ref.} \cite{Zheng17} used the bi-level structure to optimally locate charging stations to minimize travel time and energy consumption while considering traffic equilibrium.
	{Ref.} \cite{Guo18} developed a bi-level integer programming model to locate charging stations in the manner of {of the construction and deviation costs are minimized} while maximizing the number of served EV by the charging service.\\
	{Ref.} \cite{He18} proposed a bi-level programming model with the consideration of EV’s driving range, for finding the optimal locations of charging stations: the upper level {consists of optimizing} the position of charging stations to maximize the path flows that use the charging stations, while the user equilibrium of route choice with the EV’s driving range constraint is formulated in the lower level.\\
	{Ref.} \cite{Makhlouf19} developed a bilevel problem where the upper-level {objective} is a max-cover type station location and sizing problem, and the lower-level problem represents the preference of EV user behaviour in terms of making the minimum number of stops to reach their destination. \\
	{Ref.} \cite{Huang20} considered congested travel and congested stations under elastic demand to maximize profits of electric vehicle charging station owner by means of a genetic algorithm.\\
	{Ref.} \cite{Tran21} developed a bi-level program to determine the optimal location of public fast-charging stations while simultaneously considering heterogeneous vehicle classes, the installation cost of charging stations, link congestion and route choice behaviours of travellers with multiple recharging locations. 
	Recently, \cite{Piedra22} have proposed a bilevel model to combine a double public-private perspective of the problem of locating recharging stations for electric vehicles for an energy company that already has a set of supply points in the road network of a specific geographical area. At the upper level, a criterion of conditional coverage of set must be met, with which it is ensured that all users (admitting that there are no capacity restrictions in the service) can have at least one nearby refuelling point and, if that failed, they would have another alternative point at a pre-established coverage distance (thus resulting in reinforced coverage). At the lower level, once compliance with the previous requirement is guaranteed, the energy operator would be allowed to freely establish recharging points, maximizing the expected profit (according to existing demand) within the established budget. For this second criterion, the authors proposed the use of a knapsack model.\\
	{Our proposal is a more developed version of \cite{Piedra22} in which we transform the two-level model into a single-level model to compare it with a heuristic based on the knapsack and conditional coverage models.}

	\section{Model development}
	\subsection{Input data}
	In order to determine optimal location of charging stations in corridors a connected graph $G = (V,A)$ is assumed composed of a node set $V$ representing gas stations or cities and an arc set $A$ representing road sections between points, such that the existence of a shortest path in terms of distance (or travel time) between each pair of points of $V$ is always guaranteed inside $G$.
	The following notation is used in our formulation:

	Origin-destination demand matrix $(d_{ij}) , i,j \in V$.

	$(\Gamma_{ij})$: shortest path matrix between pair of nodes $i,j \in V$, where

	 $\Gamma_{ij}=\{i, v_1, v_2, \ldots,v_k,j\}$ with $v_1, v_2, \ldots,v_k$ intermediate nodes.

	$T = (t_{ij})$: distance matrix between pair of nodes through the shortest path.

	Term $q_k$ denotes the capacity of each node $k\in V$ to install charging stations,

	and $p_k$ the unit price depending on site $k$.\\
	Moreover, the following variables are required in the model:

	$x_k$: Integer variable that indicates the number of charging facilities installed

	 at point $k\in V.$

	$y_l$: Binary variable that takes value 1 if we select point $l\in V$ to open at

	 least one charging facility.\\	
	Note that the total installation cost at point $k \in V$ depends on the number of facilities $x_k$ $(0\leq x_k\leq q_k)$ that we open. It will be $x_k \cdot p_k$. Thus, the total cost in the whole network will be 
	$$\sum_{k\in V} x_k \cdot p_k$$
	\subsection{Preprocess}
	For each point $k$, the shortest paths containing point $k$ as an intermediate node are identified. This collection is labeled as $\gamma_k$.
	$$\gamma_k =\{\Gamma_{ij}; i,j \in V | k\in \Gamma_{ij}\}$$
	Once set $\gamma \equiv (\gamma_k)$ is determined, the following weights are obtained
	$$\omega_k = \sum_{\Gamma_{ij}\in\gamma_k} d_{ij}$$
	that will serve to quantify the attractiveness of locating a service at point $k$.
	\subsection{Conditional Covering and Knapsack Models }
	
	Suppose we want each favorable facility decision ($x_k \geq 1 $) to be complemented with the installation of another facility at a location $l\in V$ such that the travel time between points $k$ and $l$ $(t_{kl})$ does not exceed a quantity $R$ which has previously been set by the experts.
	\begin{definition}
		Let $G = (V,A)$ be a graph. It is called $R$-dense if $\forall p\in V, \exists q\in V (q\not= p) \quad\mbox{such that} \quad dist_G(p,g) \leq R$.
	\end{definition}
Considering that the construction costs are associated with the number of installations, the objective of the master problem is the minimization of the number of installations necessary to guarantee a reinforced coverage in the network. The Conditional Covering Problem (CCP) consists of minimizing the total number of facilities that must be established in order to cover all nodes and no facility can cover the site on which it is located, and must therefore be covered by another established facility. A site is said to be covered if its distance to the nearest facility is less than or equal to the covering radius $R$.
	
	We can observe that if $G$ is $R$-dense it guarantees us the possible existence of paths with conditional covering.\\
	The CCP was first introduced \cite{Moon84}, where an integer programming model for this problem was proposed and linear programming relaxation methods were applied to them. In \cite{Chaudhry87}, the authors consider several greedy heuristics for solving CCP and provide computational results for the same. In \cite{Moon95}, Moon and Papayanopoulos discuss a slight variation of CCP on tree graphs. In this problem, each demand point has a specific radius such that a facility has to be located within that radii. In \cite{Goldberg05} an $O(n^2)$ algorithm for the CCP on paths with a covering radius is uniform for all the vertices and arbitrary positive costs are assigned to vertices has been presented. They also improve the result with an $O(n)$ time algorithm, when the covering radius is uniform and cost is unity for all vertices of the path. In \cite{Horne05}, Horne and Smith extend the $O(n^2)$ algorithm, obtained in \cite{Goldberg05}, to the case when vertices are assigned to an arbitrary covering radius. In \cite{Benkoczi12}, new upper bounds have been proposed for the conditional covering problem on paths, cycles, extended stars, and trees. \\
	To incorporate the concept of conditional coverage to our model, we can define a parameter matrix $B = (b_{kl})$; $b_{kl} \in \{0,1\}$, such that
	$$b_{kl}=\left\{ \begin{matrix}
		1 \qquad \mbox{if} \quad t_{kl} \leq R \quad \mbox{and } k\not= l\\
		0\qquad \mbox{if} \quad t_{kl} > R \quad \mbox{or } k= l \end{matrix} \right. $$
	From this matrix $B$ we can extract the vector $B_k = \{l\in V (l\not=k) | b_{kl} = 1 \}$ for the conditional covering between pair of selected nodes.\\
	To ensure a reinforced coverage of services, the following conditional covering model is proposed
	\begin{alignat}{2}
		& \minimize
		&\qquad &  \sum_{l\in V} y_l \label{mod2}\\
		&\text{s.t.}  &  & \sum_{l\in V,l\not= k} b_{kl}  y_l \geq 1 \qquad \qquad \forall k\in V,  \label{mod2_0}\\
		&&& y_O = 1,\\
		&&& y_D = 1, \\
		&&& y_l \in \{0,1\}\qquad \qquad \quad\qquad\forall l \in V.\label{mod2_1}\\
		\nonumber
	\end{alignat}
	The inclusion of restrictions guarantees the coverage of services at the points of origin ($O$) and destination ($D$) and along any itinerary established inside the graph $G$.\\
	In addition to provide an standard level for service coverage, the conditional covering model guarantees that users requesting a nearby charging point can have a second option in the event of an incident. Thus, the CCP  can solve the problem of locating power stations from the administration point of view, since it results in a network in which all nodes are covered. They are also covered in a reinforced way to avoid cases of collapse of certain points that have high demand or cases of breakdowns or technical problems.\\
	To solve the problem from the point of view of the energy companies, the classic knapsack problem is used, which maximizes profit taking into account a maximum budget $P$.
	\begin{alignat}{2}
		& \maximize
		&\qquad &  \sum_{k\in V} \omega_k \cdot x_k \label{mod1}\\
		&\text{s.t.}  &  & \sum_{k\in V} x_k  p_k \leq P,  \label{mod1_0}\\
		&&& x_k \in \mathcal{N}, 0\leq x_k \leq q_k. \qquad \qquad \qquad\forall k \in V\label{mod1_5}\\
		\nonumber
	\end{alignat}
	The optimization model formulated corresponds to a bounded knapsack problem \cite{Kellerer04}.
	
	\subsection{Bilevel optimization model and solution algortihm}
	In practice, the goodness of the decision is based on the combination of two objectives: the minimization of the number of installations and the maximization of the attractiveness of the facilities. Taking into account the previous existence of bilevel models in the specialized literature, a bilevel programming approach with novel aspects is proposed for this problem. According to \cite{Piedra22}, the global problem can be modeled as a network optimal decision problem involving two nested objectives. \\
	In the objective function of the leading problem, the aim is to minimize the number of refuelling points necessary to guarantee a reinforced service coverage for all users. This objective is linked to the political decision maker interested in the viability of solutions.\\
	The objective function of the follower problem is to maximize the volume of demand subject to budget constraints; a criterion typically used by the private sector initiative.\\
	The full formulation is as follows:
	\begin{alignat}{2}
		&\minimize
		&\qquad &  \sum_{l\in V}  y_l \label{mod4}\\
		& \text{s.t.}  &  & \sum_{l\in V,l\not= k} b_{kl}  y_l \geq 1, \hspace*{2.2cm} \forall k\in V  \label{mod4_0}\\
		&&& \maximize  \sum_{k\in V} \omega_k \cdot x_k \label{mod4_2}\\
		&&&\text{s.t.} \quad   \sum_{k\in V} x_k  p_k \leq P,  \label{mod4_3}\\
		&&& \qquad y_k \leq x_k, \hspace*{2.7cm}\forall k \in V\label{mod4_4}\\
		&&& \qquad\frac{x_k}{q_k}\leq y_k;\quad \frac{x_k}{q_k} \leq \sum_{l\in B_k} y_l, \hspace*{0.25cm}\forall k \in V\label{mod4_5}\\
		&&& x_k \in \mathcal{N}, 0\leq x_k \leq q_k,  \hspace*{1.5cm}\forall k \in V\label{mod4_6}\\
		&&& y_l \in \{0,1\}\hspace*{3.2cm}\forall l \in V\label{mod4_7}\\
		\nonumber
	\end{alignat}
	Where (\ref{mod4}) is the leading objective of the bilevel model. 
	Constraints (\ref{mod4_0}) guarantee the reinforced coverage of all refuelling points. 
	(\ref{mod4_2}) is the follower objective of the bilevel model, adapting the number of charging points to potential demand.
	Constraints (\ref{mod4_3}) are budget constraints. 
	Constraints (\ref{mod4_4}) establish that the places selected by the leader objective have recharge points.
	Constraints (\ref{mod4_5}) ensure that all places where a recharging point is planned to be installed ($0< \frac{x_k}{q_k}\leq 1$) must be covered in a reinforced way by another different place. 
	Constraints (\ref{mod4_6})-(\ref{mod4_7}) indicate the nature of the variables used in the model.

	The Knapsack Problem is of combinatorial nature and, in computational complexity theory, is classified as NP-hard problem \cite{Garey79}. Bilevel optimization problems are known to be intrinsically hard to solve. Even the models with both linear leader and follower’s problems, which are generally the simplest to solve, are shown to be strongly NP-hard \cite{Labbe21}. Typically, solution methods used for these studies are metaheuristics, such as a genetic algorithm or large-scale neighbourhood search, or a single-level reformulation of the bilevel problems is proposed to be able to solve the models using commercial solvers.
	
	Taking these precedents into account, we propose this heuristic for solving the optimization problem in order to determine the locations and the number of charging points inside the network. \\
	\newpage
	\begin{verbatim}
		Heuristic
		1. Let S1 be the solution of the Conditional Covering Problem in V
		 [model (1),(2),(5)].
		
		2. Let S2 be the solution of the Bounded Knapsack Problem in V forcing
		that the locations that are solutions of S1 are also solutions in S2 
		[S1 is a subset of S2, obviously].
		3. Let G = {}
		4. While K = S2\S1 not empty do:
		
		4.1 Let k1 = argmax{w_i/p_i | i in K}.
		4.2 Let S1 be the solution of the Conditional Covering Problem in V,
		forcing that k1 is part of the solution [y_k1 = 1].
		4.3 Let S2 be the solution of the Bounded Knapsack Problem in V forcing
		that the locations that are solutions of S1 are also solutions in S2.
		5. End
	\end{verbatim}
	\section{Our case study}
	
\begin{figure}[h]
		\centering
		\includegraphics[scale=0.3]{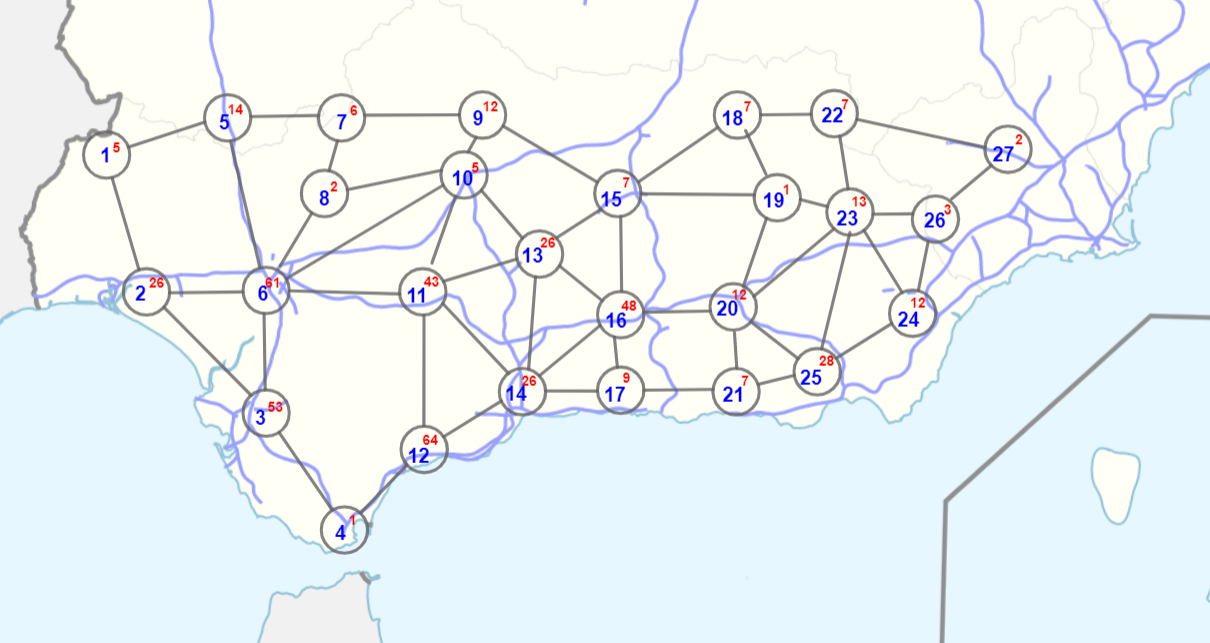}
		\caption{Example of an area with 27 existing gas stations}
		\label{Fig2}
	\end{figure}

	In order to illustrate the developed methodology, we have designed an experimental scenario inspired by a real case with existing gas stations in the southern region of Spain \cite{Repsol}, let suppose an area with 27 already installed gas stations where the distances between connections are equal to $1$, as it is shown in Figure \ref{Fig2}. Let suppose that for all sites $k$ the unit cost of each charging point are equal to 1 ($p_k = 1$), where the total budget is equal to $20$ ($P = 20$). The capacity at each node is equal to $5$ ($q_k = 5$), the attractiveness of each node is equal to the number that is shown in red in Figure \ref{Fig2}. Additionally, we assume that the solution graph is 1-dense ($R = 1$) and that one charging station is covered by another if they have a distance of 1.
	\\
	\textbf{Step 1.} If the CCP is solved with $V =\{1,2, \ldots, 27\}$, the obtained solution is $$ S1 = \{2,3,7,9,14,17,22,23\}$$
	\textbf{Step 2.} Next, the KP is solved forcing that S1 are also solutions. The solution is
	$$ S2 = \{2,3,6,7,9,12,14,17,22,23\}$$\\
	with 
	$$ x_2 = x_7 = x_{9} = x_{14}=x_{17}=x_{22}=x_{23} = 1, \quad x_3 = 3, \quad x_6 = x_{12}= 5$$

	\begin{figure}[h!]
		\centering
		\includegraphics[scale=0.5]{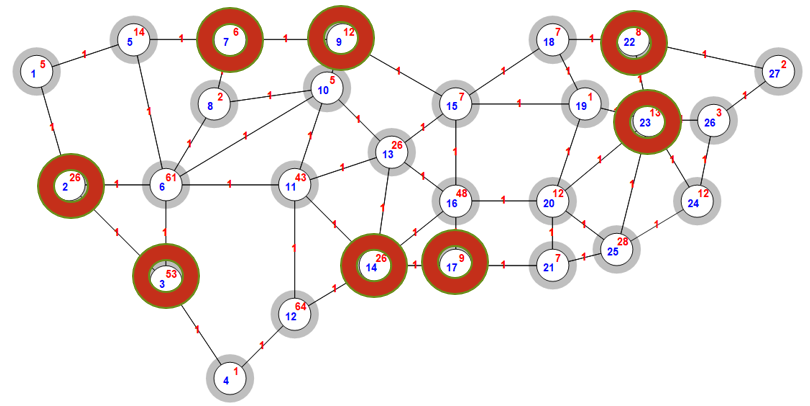}
		\includegraphics[scale=0.5]{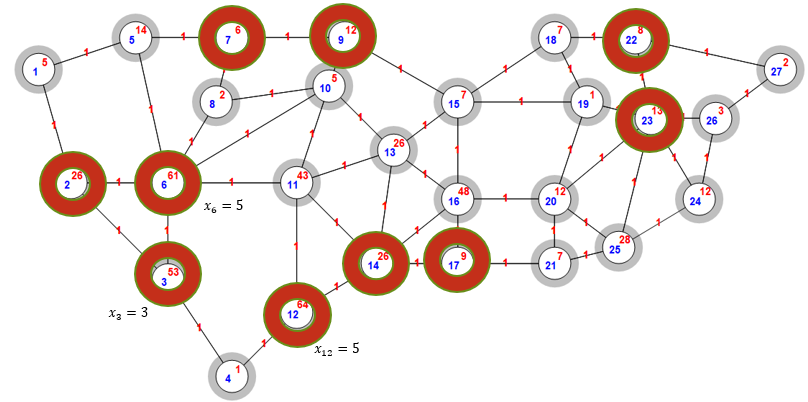}
		\caption{Selected nodes at Step 1 and Step 2.}
	\end{figure}

	\textbf{Step 3.} Once set $K = S2\setminus S1= \{6,12\}$ is determined and identified as not empty, parameter $k_1$ is assigned to node $12$. Solve the CCP including $k_1$ in the solution.
	$$S1 = \{5,6,12,13,14,15,20,22,23\}$$
	\textbf{Step 4.} Solve the KP forcing that this new version of $S1$ takes part of the solutions. 
	$$S2 = \{3,5,6,12,13,14,15,20,22,23\}$$
	\textbf{Step 5.} Repeat Step 3 for $K =S2\setminus S1= \{3\}$. Now  $k_1 := k_1 \cup \{3\}$. Solve the CCP including $k_1$ as solution.
	$$S1 = \{2,3,7,9,12,14,20,22,23\}$$
	\textbf{Step 6.} $S2 = \{2,3,6,7,9,12,14,20,22,23\}$ \\
	\textbf{Step 7.} $K = S2\setminus S1= \{6\}$ is not empty.  $k_1 = k_1 \cup \{6\}$. Solve the CCP including $k_1$ as solution.
	$$S1 = \{3,5,6,12,14,15,18,23,25,26\}$$\\
	\textbf{Step 8.} $S2 = \{3,5,6,12,14,15,18,23,25,26\}$, with 
	$$ x_5 = x_{14} = x_{15} = x_{18}=x_{23}=x_{25}=x_{26} = 1, \quad x_3 = 3, \quad x_6 = x_{12}= 5$$ \\
	\textbf{Step 9.} $K$ is empty. Stop. The Solution is $S2$.
	\begin{figure}[h!]
		\centering
		\includegraphics[scale=0.6]{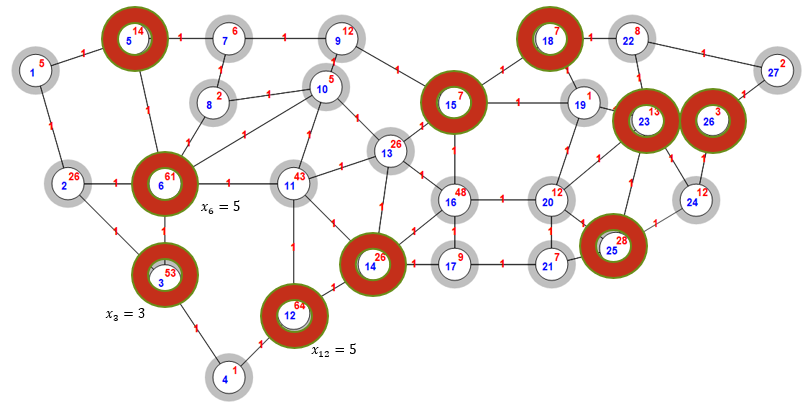}
		\caption{Selected nodes at Step 9.}
	\end{figure}\\

	Therefore, the solution obtained after applying the algorithm is to open 10 electric gas stations whose overall attractiveness would reach the value of 882. On the other hand, if we solve this same example with the exact model, the solution obtained is $S=\{2,3,7,9,14,$ $17,22,23\}$ with $x_7 = x_{9}=x_{17} = x_{22}=x_{23}=1, \quad x_2 =x_3=x_{14}=5$. In other words, the optimal solution is now to open 8 electric gas stations that give us an attractiveness of 573.
	The bilevel problem, solved by using our algorithm, has been separated into two classic linear problems and, since both are still NP-hard, we will use different heuristics from the literature to solve each one of them \cite{Lofti97,Deineko11, Pisinger00}.
	\section{Others perspectives}
	\subsection{Single level reformulation with primal-dual problems}\label{reform}
	In order to solve problem (\ref{mod4})-(\ref{mod4_7}), we reformulate it as
	a single-level optimization problem by exploiting the primal-dual
	optimality conditions of linear programming. First we define the dual problem of the follower problem knowing that the upper level binary variable $y_k$ is an input parameter whose value is fixed. 
	
	\begin{alignat}{2}
		&\minimize&\qquad &  P\cdot\lambda + \sum_{k\in V} q_k\beta_k +\sum_{k\in V} \left(\sum_{l\in B_k} y_l \right)\cdot\alpha_k+\sum_{k\in V}y_k\cdot\delta_k  + \sum_{k\in V}y_k\cdot\gamma_k \label{mod6}\\
		& \text{s.t.} && p_k \lambda + \beta_k + \frac{\alpha_k}{q_k} + \gamma_k + \frac{\delta_k}{q_k}  \geq \omega_k,\hspace*{0.75cm} \forall k \in V\label{mod6_1} \\
		&&& \lambda , \beta_k, \alpha_k, \delta_k  \geq 0,  \hspace*{3.1cm}\forall k \in V\label{mod6_2}\\
		&&& \gamma_k \leq 0. \hspace{4.6cm} \forall k \in V
	\end{alignat}
	
	According to the weak and strong duality theorems, if $x_k$ is
	a feasible solution of primal problem, $\lambda , \beta_k, \alpha_k, \delta_k,  \gamma_k$ is a feasible solution of dual problem and  the optimal values of the primal and dual problem coincide, then $x_k$ (resp. $\lambda , \beta_k, \alpha_k,\delta_k,  \gamma_k$) is an optimal solution of primal (resp. dual). Problem (\ref{mod4})-(\ref{mod4_7}) can be reformulated as the following single level problem.
	
	\begin{alignat}{2}
		{\text{E1:}} & \quad\minimize&\quad &  \sum_{l\in V}  y_l \label{mod5}\\
		& \quad\text{s.t.}  &  & \sum_{l\in V,l\not= k} b_{kl}  y_l \geq 1, \hspace*{2.8cm} \forall k\in V  \label{mod5_0}\\
		&&&   \sum_{k\in V} x_k  p_k \leq P,  \label{mod5_3}\\
		&&&  y_k \leq x_k, \hspace*{4.1cm}\forall k \in V\label{mod5_4}\\
		&&& \frac{x_k}{q_k} \leq y_k; \quad\frac{x_k}{q_k} \leq \sum_{l\in B_k} y_l, \hspace*{1.65cm}\forall k \in V\label{mod5_5}\\
		&&& p_k \lambda + \beta_k + \frac{\alpha_k}{q_k} + \gamma_k + \frac{\delta_k}{q_k}  \geq \omega_k,\hspace*{0.55cm} \forall k \in V\label{mod5_1} \\
		&&& P\lambda + \sum_{k\in V} q_k\beta_k +\sum_{k\in V} \left(\sum_{l\in B_k} y_l \right)\alpha_k+\sum_{k\in V}y_k\delta_k + \label{mod5_8}\\&&&\nonumber \hspace{5cm} + \sum_{k\in V}y_k\gamma_k =  \sum_{k\in V} \omega_k x_k,  \\
		&&& x_k \in \mathcal{N}, 0\leq x_k \leq q_k, \hspace*{2.15cm}\forall k \in V\label{mod5_6}\\
		&&& y_l \in \{0,1\},\hspace*{3.65cm}\forall l \in V\label{mod5_7}\\
		&&& \lambda , \beta_k, \alpha_k \geq 0, \hspace*{3.35cm}\forall k \in V\\
		&&& \gamma_k \leq 0. \hspace*{4.35cm} \forall k \in V
	\end{alignat}
	Constraint (\ref{mod5_8}) is the strong duality condition stating that the primal and dual objectives of the lower level problem must be equal,  { (\ref{mod5}) represents the objective function and (\ref{mod5_0}) the upper level problem constraints},  (\ref{mod5_3})-(\ref{mod5_5}) {are} the lower level primal problem constraints and (\ref{mod5_1}) represent the lower level dual problem
	constraints. 
	\subsection{Single level reformulation with exchange in the hierarchy criteria}\label{reform2}
		In view of the results of our case study in the exact model, although the solution minimizes the number of electric stations that we open, it does not give good attractiveness results. Therefore, we carry out an exchange at the hierarchy levels: 
	\begin{alignat}{2}
		&\maximize
		&\qquad & \sum_{k\in V} \omega_k \cdot x_k  \label{mod10}\\
		& \text{s.t.}  &  &\sum_{k\in V} x_k  p_k \leq P,  \label{mod10_0}\\
		&&& \frac{x_k}{q_k}\leq y_k;\quad \frac{x_k}{q_k} \leq \sum_{l\in B_k} y_l, \hspace*{0.85cm}\forall k \in V\label{mod10_5}\\
		&&& \minimize \sum_{l\in V}  y_l  \label{mod10_2}\\
		&&&\text{s.t.} \quad \sum_{l\in V,l\not= k} b_{kl}  y_l \geq 1, \hspace*{1.3cm} \forall k\in V  \label{mod10_3}\\
		&&& \qquad y_k \leq x_k, \hspace*{2.7cm}\forall k \in V\label{mod10_4}\\
		&&& x_k \in \mathcal{N}, 0\leq x_k \leq q_k,  \hspace*{1.5cm}\forall k \in V\label{mod10_6}\\
		&&& y_l \in \{0,1\}\hspace*{3.2cm}\forall l \in V\label{mod10_7}\\
		\nonumber
	\end{alignat}
Using the same argument and results used in \ref{reform}, we reformulate our problem, obtaining the following single-level model:
	\begin{alignat}{2}
	{\text{E2:}}\quad	&\maximize
	&\quad & \sum_{k\in V} \omega_k \cdot x_k  \label{mod11}\\
	& \text{s.t.}  &  &\sum_{k\in V} x_k  p_k \leq P,  \label{mod11_0}\\
	&&& \frac{x_k}{q_k}\leq y_k;\quad \frac{x_k}{q_k} \leq \sum_{l\in B_k} y_l, \hspace*{1.3cm}\forall k \in V\label{mod11_5}\\
	&&& \sum_{l\in V,l\not= k} b_{kl}  y_l \geq 1, \hspace*{2.5cm} \forall k\in V  \label{mod11_3}\\
	&&& y_k \leq x_k, \hspace*{3.8cm}\forall k \in V\label{mod11_4}\\
	&&& \sum_{l\in B_k} \alpha_k +\beta_k  \leq 1,\hspace*{2.4cm} \forall k \in V\label{mod11_1} \\
	&&&  \sum_{k\in V}\alpha_k+\sum_{k\in V}\beta_k x_k =  \sum_{k\in V} y_k \label{mod11_8},  \\
	&&& x_k \in \mathcal{N}, 0\leq x_k \leq q_k, \hspace*{2.05cm}\forall k \in V\label{mod11_6}\\
	&&& y_l \in \{0,1\},\hspace*{3.65cm}\forall l \in V\label{mod11_7}\\
	&&& \alpha_k, \beta_k \in \{0,1\}, \hspace*{3.05cm}\forall k \in V \label{mod11_8}\\
\end{alignat}
	\subsection{Computational results}
\begin{figure}[h!]
		\centering
		\includegraphics[scale=4.1]{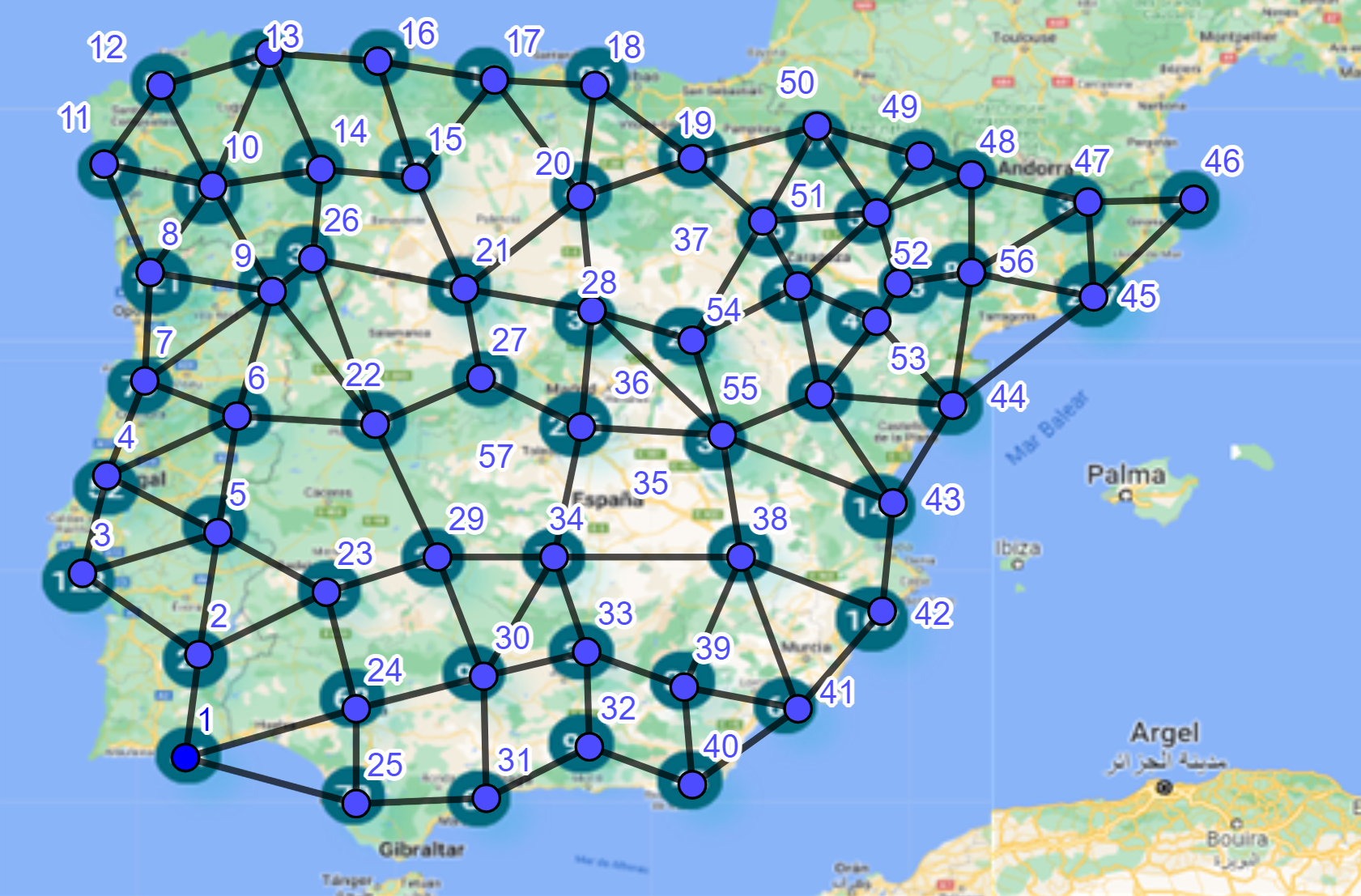}
		\caption{{Example of an area with 57 existing gas stations.}}
		\label{Fig3}
	\end{figure}

	A computational experience, consisting of performing 25 experiments on the network in Figure \ref{Fig2} has been carried out in order to compare the solutions provided in both reformulations and the heuristics solutions.
	 {Also, to check the effectiveness of the heuristic, an expansion of the network from Figure \ref{Fig2} to the network in Figure \ref{Fig3} with 57 nodes has been carried out.}
		\\
	
	 Computational experience was solved by means of the GuRoBi 9.1.2 in Python solver on a laptop with 16 GB of RAM and an Intel processor i7-1165G7 (with a 64-bit Windows 10 professional operating system and limiting the execution time to 2 h).
	
	Table \ref{Table1} shows the obtained results {by means of the three methodologies: H, E1, E2.}
\begin{table}[]
	\centering
	\begin{tabular}{|c|c|c|c||c|c|c||c|c|c|}
		\hline
		& \multicolumn{3}{c||}{H} & \multicolumn{3}{c||}{E1} & \multicolumn{3}{c|}{E2} \\ \hline
		\# & \textbf{NH }        & \textbf{AH }  &   \textbf{CPU H }      & \textbf{NE1 }        & \textbf{AE1}  & \textbf{CPU E1 }       & \textbf{NE2 }       & \textbf{AE2 } &   \textbf{CPU E2 }      \\ \hline
		1  & 10         & 862 &   \textbf{0.00}       & \textbf{8}  & 573   &  0.19     & 10         & \textbf{945}& 0.03 \\ \hline
		2  & 9          & 1481 &      \textbf{0.00}   & \textbf{8}  & 1457& 0.13        & 10         & \textbf{1627}&0.03 \\ \hline
		3  & 11         & 1600   &      \textbf{0.00} & \textbf{8}  & 747 &   0.19      & 11         & \textbf{1605}& 0.03\\ \hline
		4  & 9          & 507  &       \textbf{0.00}  & \textbf{8}  & 391 &   0.09      & 9          & \textbf{542}&  0.03\\ \hline
		5  & 10         & 849 &         \textbf{0.00}& \textbf{8}  & 573 &   0.25       & 10         & \textbf{871}& 0.02 \\ \hline
		6  & 9          & 507  &   \textbf{0.00}      & \textbf{8}  & 345 &   0.07      & 10         & \textbf{648}& 0.03 \\ \hline
		7  & 6          & 1052  &     \textbf{0.00}  & \textbf{4}  & 635 &   0.48      & 6          & \textbf{1060}& 0.09\\ \hline
		8  & 5          & \textbf{1095} &  \textbf{0.00}  & \textbf{4}  & 495 &   0.14      & 5          & \textbf{1095}&0.03 \\ \hline
		9  & 13         & 2405  &     \textbf{0.01}  & \textbf{10} & 1446 &     8.2    & 13         & \textbf{2460}& 0.08\\ \hline
		10 & \textbf{4} & 558  &    \textbf{0.02}     & \textbf{4}  & 412 &  0.78        & 5          & \textbf{673}& 0.05 \\ \hline
		11 & 8          & \textbf{1692}& \textbf{0.01}  & \textbf{4}  & 746&0.58          & 8          & \textbf{1692}&0.05 \\ \hline
		12 & -          & -  &-      & \textbf{9}  & \textbf{297} &0.07 & \textbf{9} & \textbf{297}&  \textbf{0.04}\\ \hline
		13 & 8          & 1509  & \textbf{0.01}        & \textbf{5}  & 180  &     2.99   & 8          & \textbf{1734}& 0.06 \\\hline
		14 & -          & -    & -         & \textbf{10} & 434  &   0.09     & 11         & \textbf{446}& \textbf{0.04} \\ \hline
		15 & 5          & 1248    &  \textbf{0.01}    & \textbf{4}  & 821 &      1,01   & 6          & \textbf{1580}& 0.03\\ \hline
		16 &6 &1224 &\textbf{0.00} &\textbf{4} &454 &0.35 &8 &\textbf{1287} &0.01 \\ \hline
		17 &5 &911 &\textbf{0.00} &\textbf{4 }&344 &0.07 &6 &\textbf{1033} &0.0 \\ \hline
		18 &7 &1223 & \textbf{0.00}&\textbf{5} &293 &0.62 &7 &\textbf{1223} &0.01 \\ \hline
		19 &- &- &- &\textbf{9} &7074 &0.38 &\textbf{9} &\textbf{10338} &\textbf{0.01} \\ \hline
		20 &6 & 48977&\textbf{0.00} &\textbf{4 }&38743 & 1.01&6 &\textbf{53362} &0.02 \\ \hline
		21 &6 &\textbf{15561} &\textbf{0.00} &\textbf{5} &9479 &1.64 &6 &\textbf{15561} & 0.02\\ \hline
		22 &6 &33497 &\textbf{0.00} &\textbf{4} &15471 &1.31 &6 &\textbf{40682} &0.02 \\ \hline
		23 &- &- &- & \textbf{9}&7062 &0.41 &10 &\textbf{14441} &\textbf{0.02} \\ \hline
		24 &7 &38071 &\textbf{0.00} &\textbf{4} &28073 &0,37 &7 &\textbf{39529} &0.02 \\ \hline
		25 &12 &21797 &\textbf{0.00} &\textbf{10} &7798 &0,05 &12 &\textbf{23710} &0.01 \\ \hline
		{26} &{18} &{52352} &{\textbf{0.01}} &{\textbf{16}} &{27747} &{1.96} &{19} &{\textbf{62194}
		} & {0.05} \\ \hline
		{27}& {-} &{-} &{-} &{\textbf{16}} &{9476}& {2.02}&{19} &{\textbf{24706}} & {\textbf{0.03}} \\ \hline
		{28} & {20}& {113126}&{\textbf{0.01}} &{\textbf{16}} &{88249} &{13.25} & {20}&{\textbf{114853}
		} & {0.05} \\ \hline
		{29} &{15} &{115051} &{\textbf{0.01}} &{\textbf{7}} & {29448} &{0.57} &{15} &{\textbf{116223}} & {0.02 }\\ \hline
		{30} &{22} & {112891}&{\textbf{0.01}} &{\textbf{16}} &{69394} &{13.05} & {22}&{\textbf{113754}
		} & {0.03} \\ \hline
	\end{tabular}
\caption{Results obtained from the computational experience.}
\label{Table1}
\end{table}
The columns of \textit{NH, AH} show the number of nodes and the attractivity obtained with the heuristic solution respectively. The columns of \textit{NE1, AE1} and \textit{NE2, AE2} show the number of nodes and attractivity obtained with the exact method explained in section \ref{reform} and \ref{reform2} respectively. The computation times (in seconds) of the heuristic, model 1 and model 2 are included in the columns \textit{CPU H, CPU E1, CPU E2} respectively.\\
 Instance $\#1$ is our case study. Instances between $\#2$ and $\#8$ have a single parameter randomly changed from the original example. Instances between $\#9$ and $\#18$ all the parameters have been taken randomly with values of attractiveness and budget less than 100, capacity less than 10, unit price less than 1/4 of the budget and coverage radii less than 3. Instances between $\#19$ and $\#25$ all the parameters have been taken randomly with values of attractiveness and budget less than 1000 and greater than 100, capacity less than 20 and all other parameters taken randomly in the same way. {Instances between  $\#26$ and  $\#30$ are those carried out in the network of Figure \ref{Fig3}, taking the data randomly.}\\
 In all the experiments carried out, the minimum number of nodes is given by the first exact model and the maximum attractivity is given by the second exact model. \\
The heuristic gives us a good solution between both methodologies, in all the experiments, the difference of nodes with the minimum numbers is less or equal than 4 and  with the maximum attractivity has {$5.00\%$} of relative error. The heuristic gives solutions that are closer to the optimal value of the second model because once a solution with few nodes covering the network is achieved, the heuristic increases the number of nodes to increase attractiveness.
On instance $\#12$, $\#14$, $\#19$, $\#23$ {and $\#27$} the heuristic does not give us a solution of the problem because the initial CCP solution exceeds the budget. 

{Several larger-scale computational experiments have been conducted to observe the behavior of the exact algorithms and the heuristic with instances containing a number of gas stations closer to reality (for example, only the company Repsol has more than 3000 conventional gas stations in Spain).\\	
	For E1, it is observed that it can achieve the optimal solution for instances with up to 300 points in less than 2 hours. However, for examples with more than 400 nodes, the first exact model fails to solve it in less than 2 hours and presents a gap of over 40 \%. This suggests that E1 encounters difficulties in handling larger instances due to its higher complexity associated with its dual model.\\	
	On the other hand, E2 shows impressive performance by reaching the maximum benefit solution for a network with up to 8000 nodes. It can handle larger instances more effectively than E1 due to its focus on prioritizing node attractiveness using an extension of the knapsack model.\\	
	The difference in resolution size between E1 and E2 is attributed to their distinct approaches and the characteristics of their respective models.\\	
	Finally, the heuristic H consistently provides an intermediate solution between E1 and E2 when both models can be resolved. This indicates that through its algorithm, H manages to strike a balance between coverage optimization and node attractiveness, offering a solution that represents a combination of the approaches provided by E1 and E2.}

	\section{Conclusions}
	In this work, a methodology has been developed to optimally select, among a group of candidate sites already equipped with refuelling facilities, a series of recharging points in order to guarantee that an electric vehicle can autonomously transit within a territory and, in the event of a breakdown in the selected service station, be able to count on an alternative charging station that is within a reasonable distance radius. This concept of reinforced coverage has been formulated following the conditional covering model, which minimizes the number of installations required (Criterion 1). Complementarily, the maximization of the demand that could be satisfied subject to budgetary restrictions has been a second objective, which has been formulated as an instance of the knapsack problem (Criterion 2). This second criterion has been hierarchically combined with the previous one, giving rise to a bilevel model.
	
	Since the computational complexities of both proposed models are NP-difficult, a heuristic has {finally} been designed to solve the two phases of the bilevel model separately, following an iterative scheme ({Model H}). This first solution procedure has been compared with two other methodologies based on programming a single optimization level, by integrating the primal and dual versions of the involved constraints. The first of those {exact models} has been directly obtained without modifying the hierarchy of the criteria initially considered (ie, minimizing the number of installations required as the main objective - {Model E1}), while the second {exact model} pursues the maximization of the user population covered as main objective ({Model E2}).\\
	{E1, being a model that prioritizes coverage optimization, yields good results in terms of a smaller number of nodes, which aligns with the viewpoint of the government administration that must assist in funding new infrastructures. On the other hand, E2 prioritizes the attractiveness of the nodes using an extension of the knapsack model. Generally, increasing the number of nodes beyond the minimum required yields good results in achieving the established objective and efficient computation in obtaining solutions. The results of heuristic H have been in an intermediate position between the optimal solutions obtained through the previous exact models, determining a reasonably low number of nodes required to provide an attractive coverage for the new service.}
	
	\section*{Acknowledgements}
	This work was in part supported by the Proyectos I+D+i FEDER Andalucía 2014-2020 under grant US-1381656, and by the Ministerio de Investigación (Spain)/FEDER under grant PID2019-106205 GB-I00. This support is gratefully acknowledged.

\end{document}